\documentclass[11pt]{article}

\usepackage{amsmath}
\usepackage{amsfonts}
\usepackage{amssymb}
\usepackage{amsthm}
\usepackage{enumitem}
\usepackage{mathtools}
\usepackage{setspace}
\usepackage{etoolbox}
\usepackage{changepage}

\newtheoremstyle{mystyle}
{11pt}				
{11pt}				
{}					
{}					
{\bfseries}			
{}					
{5.5pt}				
{}					

\theoremstyle{mystyle}
\newtheorem{theorem}{Theorem}[section]

\newtheorem{proposition}[theorem]{Proposition}
\newtheorem{corollary}[theorem]{Corollary}

\renewenvironment{proof}[1][Proof.]{\vspace{-16.5pt} \begin{trivlist}
	\item[\hskip \labelsep {\bfseries #1}]}{\qed \end{trivlist}}

\setitemize{noitemsep, topsep=5.5pt, parsep=5.5pt, partopsep=0pt}

\allowdisplaybreaks
\appto\normalsize{
	\abovedisplayskip=5.5pt plus 2pt minus 2pt
	\belowdisplayskip=5.5pt plus 2pt minus 2pt
	\abovedisplayshortskip=5.5pt plus 2pt minus 2pt
	\belowdisplayshortskip=5.5pt plus 2pt minus 2pt}
\appto\small{
	\abovedisplayskip=5.5pt plus 2pt minus 2pt
	\belowdisplayskip=5.5pt plus 2pt minus 2pt
	\abovedisplayshortskip=5.5pt plus 2pt minus 2pt
	\belowdisplayshortskip=5.5pt plus 2pt minus 2pt}

\newcommand{\gap}{\vspace{11pt}}
\newcommand{\T}{\mathsf{T}}
\newcommand{\R}{\mathbb{R}}
\newcommand{\Rn}{\mathbb{R}^{n}}
\newcommand{\Rnp}{\mathbb{R}_{+}^{n}}

\newcommand{\tr}{\operatorname{tr}}

\newcommand{\spn}{\operatorname{span}}
\newcommand{\ext}{\operatorname{ext}}

\newcommand{\cone}{\operatorname{cone}}

\newcommand{\one}{\mathbf{1}}

\newcommand{\ip}[2]{\left< #1,\, #2 \right>}
\newcommand{\abs}[1]{\left\vert #1 \right\vert}

\newcommand{\opnorm}[1]{\left\vert\kern-0.25ex\left\vert\kern-0.25ex\left\vert #1 
	\right\vert\kern-0.25ex\right\vert\kern-0.25ex\right\vert}

\newcommand{\set}[2]{\left\{ #1 \ \middle| \ #2 \right\}}

\title{\bf Permutation invariant proper polyhedral cones \\and their Lyapunov rank}
\author{
	Juyoung Jeong\\
	Department of Mathematics and Statistics\\
	University of Maryland, Baltimore County\\
	Baltimore, Maryland 21250, USA\\
	juyoung1@umbc.edu\\[11pt]
	and\\[11pt]
	M. Seetharama Gowda\\
	Department of Mathematics and Statistics\\
	University of Maryland, Baltimore County\\
	Baltimore, Maryland 21250, USA\\
	gowda@umbc.edu
}

\date{\today}

\begin{document}

\maketitle

\begin{abstract}
	The Lyapunov rank  of a proper cone $K$ in a finite dimensional real 
Hilbert space 
is defined as the dimension of the space of all 
Lyapunov-like transformations on $K$, or equivalently, the dimension of the Lie algebra of the automorphism 
group of $K$. This (rank) measures the number of linearly independent bilinear relations needed to 
express a complementarity system on $K$ (that arises, for example, from a linear program or a complementarity problem on the cone). 
Motivated by the problem of describing  spectral/proper cones where the complementarity system can be expressed as a square system (that is, where  
the  Lyapunov rank is greater than equal to the dimension of the ambient space), we consider proper polyhedral cones in $\Rn$ that are permutation invariant. For such cones we show that    
the Lyapunov rank is either $1$ (in which case, the cone is irreducible) or $n$ (in which case, the cone is isomorphic to $\Rnp$). 
In the latter case, we show that the corresponding  spectral cone is isomorphic to a symmetric cone.
 \end{abstract}

\vspace{1cm}

\gap


\section{Introduction}
In optimization, from the perspective of conic programming, 
the nonnegative orthant, the semidefinite cone, and more generally, symmetric cones are considered `good' cones as 
they admit polynomial-time algorithms. These are also considered `good' in complementarity theory, as 
any corresponding complementarity system/problem could be written as a square system (thus allowing Newton type or other iterative schemes). The ability to write a complementarity system on a cone as a square system could be quantified in terms of the so-called Lyapunov rank and one could define `good' cones 
as those where the Lyapunov rank  is greater than or equal to the dimension of the ambient space.
Our motivation comes from the problem of describing such cones.

Given a Euclidean Jordan algebra $V$ of rank $n$ \cite{faraut-koranyi}, the symmetric cone $V_+$ of  $V$ is just the cone of squares in $V$. 
It consists precisely those elements of $V$ with every eigenvalue nonnegative. That is, 
$V_+=\lambda^{-1}(\Rnp),$
where $\lambda:V\rightarrow \Rn$ denotes the eigenvalue map that takes any element $x$ in $V$ to its eigenvalue vector whose entries are the eigenvalues of $x$ written in the decreasing order. 
Thus, $V_+$ is an example of a  `spectral cone' obtained as the inverse image (under the eigenvalue map) of a permutation 
invariant proper cone in $\Rn$ \cite{jeong-gowda-spectral cone}, \cite{jeong-gowda-spectral set}. 
Can one replace $\Rnp$ by other permutation
invariant proper (polyhedral) cones in $\Rn$ and generate `good' cones? This question leads to the problem of describing Lyapunov ranks of spectral cones in  Euclidean Jordan algebras, especially those induced by permutation invariant proper polyhedral cones. 
In the (basic) setting of $\Rn$, this reduces to finding the Lyapunov rank of a  permutation invariant proper polyhedral cone.  
This is the main focus of the paper.
\\
The {\it Lyapunov rank of a proper cone} $K$ in a real finite dimensional Hilbert space $V$ is defined as follows.  
Let $K^*=\{y\in V: \langle y,x\rangle \geq 0,\,\,\forall\,x\in K\}$ denote the dual of $K$. 
 A linear transformation $L : V \rightarrow V$ is said to be a {\it Lyapunov-like transformation} on $K$ if 
$$ x \in K, \ s \in K^*, \ \ip{x}{s} = 0 \Longrightarrow \ip{L(x)}{s} = 0. $$
The set of all such transformations on $K$ is a  real vector space, denoted by  LL$(K)$. 
We define the {\it Lyapunov rank of $K$} by
$$\beta(K):=\dim\,\mbox{LL}(K).$$ 
It is known that $\beta(K)$ is the dimension of the  Lie algebra of the automorphism group of $K$ (consisting of invertible 
linear transformations keeping the cone invariant) \cite{gowda-tao}.
To see why this number is useful, let us assume without loss of generality that $V=\Rn$ and consider 
the  complementarity problem CP$(f,K)$ corresponding to a mapping $f :\Rn  \rightarrow \Rn$:   
find $x,s \in \Rn$ such that
\begin{equation} \label{eq1}
	x \in K, \ s = f(x) \in K^*, \ \ip{x}{s} = 0.
\end{equation}
Such a problem appears, for example, in the study of equilibrium  and variational inequality problems, and arises in numerous applications in engineering, sciences, and economics \cite{facchinei-pang}. 
There are various strategies for solving complementarity problems, see \cite{facchinei-pang}. In the above 
complementarity system, there are $2n$ real variables $x_i, s_i$ ($i=1,2,\ldots,n$), while there are $n+1$ equations, namely, $s = f(x)$ and $\ip{x}{s} = 0$. Thus, in order to see if the above system could be written as 
a square system (which will help us to analyze the system or describe an algorithm to solve it), we rewrite the last bilinear relation $\ip{x}{s} = 0$ as an equivalent system of $\beta(K)$  linearly independent bilinear relations and check if $\beta(K)\geq n$.
Although this latter property does not always hold, it will be useful to identify cones where this is possible. 
\\
The concept of Lyapunov rank was introduced in \cite{rudolf et al} using the term `bilinearity rank'. Because these `bilinearity relations' are just Lyapunov-like transformations \cite{gowda-tao}, we have used the term `Lyapunov rank'.
The  concept of Lyapunov-like transformation 
was first introduced in \cite{gowda-sznajder} as a generalization of the Lyapunov transformation 
$X \mapsto AX + XA^{\mathrm{T}}$ that appears in the linear dynamical system theory, and has been the subject of 
several recent works. The results and properties related to  Lyapunov-like transformations can be found 
in \cite{gowda-sznajder}, \cite{gowda-tao}, and \cite{gowda-trott}.
\\
Beginning with \cite{rudolf et al}, there are a number of works related to Lyapunov rank, see the references. 
For the record, we  list below some known results \cite{gowda-tao}, \cite{gowda-trott}, \cite{orlitzky-gowda}:

\begin{itemize}
	\item[(1)] For any proper cone $K$ in  $\Rn$, $\beta(K) \leq (n-1)^2$.
	\item[(2)] For any proper polyhedral cone in $K$ in $\Rn$, $1 \leq \beta(K) \leq n$, $\beta(K) \neq n-1$.
	\item[(3)] In $\mathcal{S}^n$, $\beta(\mathcal{S}_+^n) = n^2$ and $\beta(\mathcal{CP}_n) = n$ 
where $\mathcal{S}^n_+$ and $\mathcal{CP}_n$ denote the set of all $n \times n$ positive semidefinite and completely positive matrices, respectively. 
	\item[(4)] In $\Rn$ ($n \geq 3$), for any $p\in [1,\infty]$, let $l_{p,+}^n := \lbrace x = (x_0, \bar{x}) \in \R \times \mathcal{R}^{n-1} : x_0 \geq ||x||_p \rbrace$. Then, $\beta(l_{p,+}^n) = 1$ for $p\neq 2$ and $\beta(l_{2,+}^n) = (n^2 - n + 2)/2$.
\end{itemize} 

In this paper, we raise the issue of describing the Lyapunov rank of a  proper polyhedral  cone $Q$ in $\Rn$ that is {\it permutation invariant} which means that 
$\sigma(Q)=Q$ forall $\sigma\in \Sigma_n,$
where $\Sigma_n$ denotes the set of all permutation matrices on $\Rn$.
Strengthening Item (2) above, we  show that for such a cone, $\beta(Q)$ is either $1$ (in which case, the cone is irreducible) or $n$ (in which case, the cone is isomorphic to $\Rnp$). 
When a 
permutation invariant proper polyhedral cone is isomorphic to $\Rnp$, we show that the spectral cone  
$\lambda^{-1}(Q)$ in a Euclidean Jordan algebra $V$ of rank $n$ 
is isomorphic to $V_+$
thus yielding a `good' cone (albeit, a known one).
 The problem whether such a statement holds for some/all permutation invariant proper polyhedral cones in $\Rn$ with $\beta(Q)=1$ and the
 more general problem of describing Lyapunov rank of a  permutation invariant proper cone  are left for further study. 
\section{Preliminaries}

Throughout this paper, $V$ denotes a finite dimensional real Hilbert space. 
For a set $S$ in $V$, $\spn(S)$ denotes the span of $S$ (which is the subspace generated by $S$); we denote the interior of $S$ by $\mathrm{int}(S)$. 

A nonempty set $K$ in $V$ is a {\it convex cone} if $tx+sy \in K$ for all $x,y \in K$ and $t,s \geq 0$. (So, we require every convex cone to contain zero.) 
	A convex cone $K$ is said to be {\it pointed} if $K \cap (-K)=\{0\}$ and {\it solid} if $\mathrm{int}(K) \neq \emptyset$. A {\it proper cone} in $V$ is a  closed convex cone that is pointed and solid.
We say that two convex cones are isomorphic if one can be mapped onto the other by an invertible linear transformation.

Given a nonempty set $S$ in $V$, the {\it conic hull} of $S$ is
\[ \cone(S) = \set{\sum_{i=1}^{k} \alpha_ix_i}{k \in \mathbb{N},\,x_i \in S,\, \alpha_i \geq 0}. \]

	A closed convex cone $K$ is a {\it polyhedral cone} if it is finitely generated, that is, $K=\cone(S)$ for some finite set $S$.

Let $K$ be a convex cone in $V$. A nonzero vector $x \in K$ is called an {\it extreme vector} of $K$ if $x = y+z$, with $y,\, z \in K$, implies that $y$ and $z$ are both nonnegative scalar multiples of $x$. We say that two extreme vectors  are equivalent (and hence consider them to be the `same') if they are positive scalar multiples of each other. Let 
\[ \ext(K): = \set{x \in K}{\text{$x$ is an extreme vector of $K$}}. \]

Thanks to a well-known theorem of Minkowski, for any proper cone $K$, $\ext(K)$ is nonempty and $K=\cone(\ext(K)).$

A  convex cone $K$ in $V$ is said to be  {\it reducible} if there exist nonempty nonzero sets $K_1,\, K_2 \subseteq V$ such that
\[ K = K_1 + K_2, \quad \mbox{and}\quad \spn(K_1) \cap \spn(K_2) = \{0\}. \] 
(The sets $K_1$ and $K_2$ turn out to be  convex cones.) When this happens, we say that $K$ is the {\it direct sum} of $K_1$ and $K_2$. If $K$ is proper, then  $K_1$ and $K_2$ are proper 
in their respective spans; hence, we can define their Lyapunov ranks.  
We say that a nonzero convex cone $K$ is {\it irreducible} if it is not reducible. We recall the following result. 

\begin{proposition}(\cite{hauser-guler}, Theorem 4.3) \label{decomposition of  reducible cone} {\it Any nonzero 
reducible pointed convex cone $K$ can be written as 
$$ K = K_1 + K_2 + \cdots + K_r,$$ 
where each $K_i$ is a nonzero irreducible pointed convex cone, and $\spn(M) \cap \spn(N) = \{0\}$ whenever $M$ is the  sum of some $K_i$s and $N$ is the  sum of the rest. Moreover, these $K_i$s are unique and the above representation is unique 
up to permutation of indices.}
\end{proposition}

In the setting of the above result, we  say that $K$ is the direct sum of $K_i$ and write $K=
 K_1 \oplus K_2 \oplus \cdots \oplus K_r.$
We then  have
\begin{itemize}
\item [$\bullet$]
$\ext(K)=\bigcup_{i=1}^{r}\,\ext(K_i)$ and 
\item [$\bullet$] 
when $K$ is proper, $\beta(K) = \beta(K_1) + \beta(K_2) + \cdots + \beta(K_r).$ 
\end{itemize}

\gap

An $n \times n$ {\it permutation matrix} is a matrix obtained by permuting the rows of an $n \times n$ identity matrix. The set of all $n \times n$ permutation matrices is denoted by $\Sigma_n$. For convenience, we treat 
an element $\sigma \in \Sigma_n$ either as a permutation matrix or as 
a permutation of indices $\{1,\, 2,\, \ldots,\, n\}$.
Recall that a convex cone $Q$ in $\R^n$ is {\it permutation invariant} if $\sigma(Q) = Q$ for all $\sigma \in \Sigma_n$. It is easy to generate such cones:
for any nonempty  subset $S$ in $\Rn$, 
$\cone ( \Sigma_n(S))$ is a permutation invariant convex cone.  
Some important permutation invariant  polyhedral cones are \cite{jeong-gowda-spectral cone}: $\Rnp$ and  
$$Q_{p}^n:=\{u\in \Rn:\,u_n^{\downarrow}+u_{n-1}^{\downarrow}+\cdots+u_p^{\downarrow}\geq 0\},$$
for any fixed natural number $p$, $1\leq p\leq n$.
Here for any vector $u\in \Rn$, $u^\downarrow$ denotes the decreasing rearrangement of $u$.

\section{A characterization of  permutation invariant pointed reducible cones}
Before we present our Lyapunov rank results, we characterize permutation invariant pointed reducible cones in $\Rn$.
The following result may be of independent interest.
In this result, we let $I$ denote the identity matrix and $E$ denote the matrix of all ones.

\begin{theorem} \label{theorem 1}
{\it 	Suppose $Q$ is a pointed convex cone in $\Rn$ that is permutation invariant and  reducible.
Then, $Q$ is isomorphic to $\Rnp$. 
In fact, $Q=A(\Rnp)$, where 
$A=(a-b)I+bE$, with $a,b\in \R$,  $a\neq b$, $a\neq (1-n)b$.
}\end{theorem}

\begin{proof}
As $Q$ is  reducible, we can write 
        \[ Q = Q_1 \oplus Q_2 \oplus \cdots \oplus Q_r, \]
        where $r>1$ and each $Q_i$ is a nonzero irreducible pointed convex cone.
By Proposition \ref{decomposition of  reducible cone}, these $Q_i$s are unique and the above representation is unique up to permutation of indices.	
Now, by permutation invariance of $Q$, for any 	
$\sigma\in \Sigma_n$,
$$Q=\sigma(Q)=\sigma(Q_1)\oplus\sigma(Q_2)\oplus\cdots\oplus\sigma(Q_r)$$ is another decomposition of $Q$. Hence, for each $i$, there is a (unique) $j$ such that $\sigma(Q_i)=Q_j$.
Thus, every $Q_i$ is (permutation) isomorphic to some $Q_j$. This sets up a binary relation between  indices in $\{1,2,\ldots, n\}$. Since $\Sigma_n$ is a group, this relation becomes an equivalence relation. Partitioning this index set or 
grouping isomorphic $Q_i$s together, we may write
	\[ Q = E_1 \oplus E_2 \oplus \cdots \oplus E_s, \] 
	where each $E_i$ is a direct sum of (permutation) isomorphic cones. We now claim that each $E_i$ is permutation invariant. Without loss of generality, consider  
	\begin{equation} \label{eq5}
	E_1 = Q_1 \oplus Q_2 \oplus \cdots \oplus Q_{r_1}.
	\end{equation} 
	As $\{ \sigma(Q_1),\, \ldots,\, \sigma(Q_{r_1}) \} \subseteq \{ Q_1,\, \ldots\, Q_{r_1} \}$ for any $\sigma \in \Sigma_n$, we see that $\sigma(E_1) \subseteq E_1$ for any $\sigma \in \Sigma_n$.
Thus, $E_1$, and more generally any $E_i$, is permutation invariant.	
	Now Lemma 9.2 in \cite{jeong-gowda-spectral cone} says that any (nonempty) nonzero pointed permutation invariant convex cone contains either $\one$ (the vector of ones in $\Rn$) or $-\one$. Because $Q$ and all $E_i$ are nonempty, nonzero, pointed permutation invariant cones, {\it we may assume that 
$\one$ belongs to all of them.} Since zero is the only common element in all $E_i$s, we must have $s=1$. 
Thus, $Q = E_1 = Q_1 \oplus \cdots \oplus Q_{r_1}$, where all the $Q_i$s are (permutation) isomorphic.
To simplify the notation, let $Q = Q_1 \oplus Q_2 \oplus \cdots \oplus Q_r$, where all the $Q_i$s are (permutation) isomorphic.
An immediate consequence is the following: As $\dim(Q_i) = \dim(Q_j)$ for every $i,j \in \lbrace 1,2,\ldots,r \rbrace$,
$n\geq \dim(Q) = \sum_{i=1}^{r}\dim(Q_i) = r \dim(Q_1)$
and so
$r$ divides $\dim(Q).$ 
\\
\noindent{\it Claim:} $r=n$.\\
Suppose, if possible, $1<r<n$. 
We consider the following cases and in each case, we derive a contradiction.
\\
The case $n=2$ is not possible, as $1<r<n$.\\
{\it Suppose $n=3$.} Then $r=2$ and $\dim(Q)$ (which is is divisible by $r$ and less than or equal to $n$) must be $2$. 
In this case, $Q$ is a direct sum of one dimensional (pointed) isomorphic  cones $Q_1$ and $Q_2$.    
Then $Q$ will have exactly two extreme vectors. Since one of these must be a 
nonconstant vector (and consequently will have at least two distinct entries), permuting the 
entries of this vector will result in at least three distinct extreme vectors of $Q$, leading to a contraction. (Here, we have used the fact that any permutation $\sigma$ takes an extreme vector to an extreme vector.) 
\\{\it Let $n=4$ so that  $r$ is $2$ or $3$.} \\
Suppose $r=2$. Then, as $r$ divides $\dim(Q)$, $Q$ must have dimension $2$ or $4$. 
If $\dim(Q)=2$, then $Q$ is a direct sum of one dimensional (pointed) isomorphic  cones $Q_1$ and $Q_2$.
We argue as in the previous case: $Q$ will have exactly two extreme vectors, one of which is a nonconstant vector.
Permuting the entries of this vector will result in at least three distinct extreme vectors of $Q$, leading to a contraction.
\\
If $\dim(Q)=4$, then 
$Q$ is a direct sum of two $2$-dimensional (pointed) isomorphic  irreducible cones $Q_1$ and $Q_2$. This is not possible, as there are no irreducible pointed cones of dimension 2.
\\
We now suppose $r=3$. In this case, $\dim(Q)=3$ and so, $Q$ is a  direct sum of three, one- dimensional (pointed) isomorphic  cones. Then $Q$ will have exactly three extreme vectors. Since one of these must be a
nonconstant vector, permuting the entries of this vector will result in at least four  distinct extreme vectors of $Q$, leading to a contraction.
\\Thus we have shown that the assumption $1<r<n$ is not possible when $n$ is $2$, $3$ or $4$.
\\
Now {\it suppose  $n \geq 5$.}
 We know that, for each $\sigma \in \Sigma_n$ and $Q_i$, we have $\sigma(Q_i) = Q_j$ for some unique $Q_j$. Define $\phi : \Sigma_n \to \Sigma_r$ which takes $\sigma \in \Sigma_n$ to $\phi(\sigma) \in \Sigma_r$ such that
	\[ \sigma(Q_i) = Q_{\phi(\sigma)[i]} \forall\; i=1,\,2,\, \ldots,\, r, \]
	where $\phi(\sigma)[i]$ denotes the image of index $i$ under the permutation $\phi(\sigma)$. We now show that $\phi$ is a group homomorphism. Let $\sigma_1,\, \sigma_2 \in \Sigma_n$. Then we have $\sigma_1\sigma_2(Q_i) = Q_{\phi(\sigma_1\sigma_2)[i]}$ as well as
	\[ \sigma_1\sigma_2(Q_i) = \sigma_1(Q_{\phi(\sigma_2)[i]}) = Q_{\phi(\sigma_1)[\phi(\sigma_2)[i]]} = Q_{(\phi(\sigma_1)\phi(\sigma_2))[i]}, \]
	for all $i = 1,\, 2,\, \ldots,\, r$. This shows that $\phi(\sigma_1 \sigma_2) = \phi(\sigma_1) \phi(\sigma_2)$ implying that $\phi$ is a group homomorphism. Hence, $\ker(\phi)$ is a normal subgroup of $\Sigma_n$. However, as $n \geq 5$, 
\begin{center}
{\it $\Sigma_n$ has only three normal subgroups, namely, $\{I\}$, $\Sigma_n$, and $A_n$,}
\end{center}
where $I$ represents the identity permutation/matrix and 
$A_n$ is the alternating group (of all even permutations in $\Sigma_n$), see  2.T.5 in \cite{dixon}.
\\
We now  consider the following   three possibilities:
	
	\textbf{Case 1.} Suppose $\ker(\phi) = \{I\}$. Then $\phi$ is injective and hence comparing cardinalities, we get  $\abs{\Sigma_n} \leq \abs{\Sigma_r}$. This implies $n \leq r$. Since $r<n$, this is not possible.
\\	\textbf{Case 2.} Suppose $\ker(\phi) = \Sigma_n$. In this case, we have $\phi(\sigma)[i] = i$ for all $i$. This means that for each $i = 1,\,2,\, \ldots,\, n$, we have
	\[ \sigma(Q_i) = Q_i \quad \mbox{forall}\,\, \sigma \in \Sigma_n. \]
	However, this is not possible, as all $Q_i$s are (permutation) isomorphic. \\
	\textbf{Case 3.} Suppose $\ker(\phi) = A_n$. Then $\phi(\sigma)[i] = i$, or equivalently,
$\sigma(Q_i)=Q_i$  for any even permutation $\sigma \in A_n$. Now, $\one$ (that belongs to  $Q$) can be decomposed as 
	\[ \one = u_1 + u_2 + \cdots + u_r, \]
	where $u_i \in Q_i$ for all $i$; we assume without loss of generality that 
$u_1 \neq 0$. Let 
$u_1 = (u_{11},\, u_{12},\, \ldots,\, u_{1n})^T$, where $n\geq 5$.  For the even permutation
	\[ \rho = \begin{pmatrix}
	1 & 2 & 3 & 4 & \cdots & n \\ 2 & 3 & 1 & 4 & \cdots & n
	\end{pmatrix} \]
	we have
	\[ \one = \rho(\one) = \rho(u_1) + \rho(u_2) + \cdots + \rho(u_r). \]
	Notice that for all $i$, $\rho(u_i) \in \rho(Q_i)=Q_i$. Hence, by the uniqueness of a decomposition, we have $\rho(u_i) = u_i$ for each $i$, especially $\rho(u_1) = u_1$. This implies that $u_{11} = u_{12} = u_{13}$. Similarly, by using all even permutations which permute only three entries of $u_1$, we  verify that all entries of $u_1$ are the same. This means that $u_1$ is a nonzero multiple of $\one$ and consequently (as $Q$ is pointed)  $\one \in Q_1$. However, as all $Q_i$ are (permutation) isomorphic, we must have $\one \in Q_i$ for all $i = 1,\,2,\, \ldots,\, r$, leading to a contradiction.
	
	Thus, we have shown that $1<r<n$ is not possible, proving the claim that $r=n$. This shows that each $Q_i$ has dimension one and that $Q$ is isomorphic to $\Rnp$.\\ 
Now let  $Q=B(\Rnp)$ for some  invertible $B\in \R^{n\times n}$. We  show that columns of $B$ can be permuted to get a matrix of the form $A:=(a-b)I+bE$ where $a,b\in \R$ with $a\neq b$, $a\neq (1-n)b$. This $A$ will then satisfy $Q=A(\Rnp)$.  
\\
Since $\ext(Q)$ is the set of all columns of $B$ and $Q$ is permutation invariant, permuting the entries of  any column of $B$ will result in  another column of $B$. As $B$ is invertible, $B$ cannot just have all constant columns. Among all nonconstant columns of $B$, let $a$ be a real number that 
appears the least number of times. Without loss of generality, let this column be of the form
$$u:=[a,a,\ldots, a, *,*,\ldots, *]^T,$$
where $a$ appears $m$-times ($1\leq m<n)$ and $*$ denotes {\it any} number different from $a$. \\
\noindent{\it Claim:} $m=1$ and all $*$s are equal.\\
As this statement is obvious for $n=2$, we assume that $n>2$.\\
Suppose, if possible, $m>1$. 
Let $\sigma(i,j)$ denote a permutation/transposition which interchanges the 
indices $i$ and $j$ while keeping all others unchanged. Applying permutations $\sigma(1,m+1), \sigma(1,m+2),\ldots, \sigma(1,n)$ to $u$ will result in $n-m$ distinct columns of $B$ each different from $u$.  Applying permutations $\sigma(2,m+1), \sigma(2,m+2),\ldots, \sigma(2,n)$ to $u$ will result in $n-m$ distinct columns of
$B$ that are different from the previously generated columns. More generally, applying permutations $\sigma(k,l)$ to $u$, where
$1\leq k\leq m$ and $m+1\leq l\leq n$, we can  generate distinct columns of $B$ totaling $m(n-m)$ columns all different from $u$. This leads to the inequality $1+m(n-m)\leq n$. As $m>1$, a simple calculation gives $n=m+1$. 
But then, $a$ appears $n-1$ times and $*$ appears once. As $n>2$, we reach a contradiction to the choice of $a$.
Hence, $m=1$. Now, $u$ is of the form
$$u=[a,*,*,\ldots, *]^T.$$ 
Then, applying permutations $\sigma(1,2), \sigma(1,3),\ldots,\sigma(1,n-1)$ to $u$, we get $n-1$ distinct columns of $B$ which together with $u$ yield all columns of $B$. If two $*$s are different, interchanging these will result in a column  that does not appear in the previous listings. Hence, all $*$s are equal. So, we may  write
$$u=[a,b,b,\ldots, b]^T$$
where $b\neq a$. Now, let $A=(a-b)I+bE$. Clearly, $A$ is obtained by permuting the columns of $B$. Hence, $Q=B(\Rnp)=A(\Rnp)$ and $A$ is invertible.  
 We have already observed 
that $a\neq b$. Note that $b=0$ is an acceptable value for $A$. 
When $b\neq 0$, $A=b\big (E+\frac{a-b}{b}I).$ Since the eigenvalues of $E$ are $0$ (with multiplicity $n-1$) and $n$, the eigenvalues of $E+\frac{a-b}{b}I$ are $\frac{a-b}{b}$ (with multiplicity $n-1$) and $n+\frac{a-b}{b}$. As $A$ is 
 invertible, we must have $n+\frac{a-b}{b}\neq 0$, that is, $a\neq b(1-n)$. Hence,  
$A=(a-b)I+bE$, where $a,b\in R$ with $a\neq b$, $a\neq (1-n)b$.
\\We end the proof by observing that if $A=(a-b)I+bE$ with $a\neq b$ and  $a\neq (1-n)b$, then
$Q=A(\Rnp)$ is a permutation invariant proper polyhedral cone in $\Rn$. 
\end{proof}

\gap

\noindent{\bf Remarks.} The above proof reveals the following: If a real square  matrix $B$ is invertible and 
the cone $Q=B(\Rnp)$ is permutation invariant, then there is a permutation $\sigma\in \Sigma_n$ such that 
$$Bx=(a-b)\sigma(x)+b\,tr(x)\one\quad (x\in \Rn).$$
It is interesting to observe that the above $B$ is a particular instance of an `isotone' linear map on $\Rn$ (which is a `majorization preserving' linear map), see \cite{ando}, Corollary 2.7.

\section{Lyapunov rank of a permutation invariant proper polyhedral cone}

\begin{theorem}\label{theorem 2}
{\it Suppose $Q$ is a proper polyhedral cone in $\Rn$ that is permutation invariant.
Then,
either 
\begin{itemize}
\item [(a)] $Q$ is irreducible, in which case,  $\beta(Q)=1$, or
\item [(b)] $Q$ is reducible, in which case, $Q$ is isomorphic to $\Rnp$ and
$\beta(Q)=n$.
\end{itemize}
}
\end{theorem}

\begin{proof}
As $Q$ is a proper polyhedral cone, $(a)$ follows from \cite{gowda-tao}, Corollary 5.\\
Now suppose $Q$ is reducible. As $Q$ is permutation invariant and pointed, from Theorem \ref{theorem 1},
$Q$ is isomorphic to $\Rnp$; hence $\beta(Q)=\beta(\Rnp)=n$. 
\end{proof} 
 
The following are easy to verify.

\begin{corollary}
{\it Let $Q$ be a  permutation invariant proper polyhedral  cone in $\R^n$. If $Q$ has more than $n$ extreme vectors, then $Q$ is irreducible and  $\beta(Q) = 1$.
}\end{corollary}

\begin{corollary}
{\it        Let $Q$ be a  permutation invariant proper polyhedral   cone in $\R^n$. Suppose there exist $d_1,\, d_2 \in \ext(Q)$ such that $\sigma(d_1) \neq d_2$ for all $\sigma \in \Sigma_n$. Then $Q$ is irreducible and  $\beta(Q) = 1$.
}
\end{corollary}

\section{Spectral cones induced by  permutation invariant proper polyhedral cones}
Let $V$ be an Euclidean Jordan algebra of rank $n$  with $\lambda:V\rightarrow \Rn$ denoting the eigenvalue map (see the Introduction). For any permutation invariant proper polyhedral cone $Q$ in $\Rn$,  consider the {\it spectral cone} defined by 
$$K:=\lambda^{-1}(Q).$$

In the result below, we partially address the question of when the above spectral cone is `good'.
Here, $e$ denotes the unit element in $V$; for any $x\in V$, $tr(x)$ denotes the trace of $x$ which is the sum of all eigenvalues of $x$.

\begin{proposition}
{\it Let $V$ be any Euclidean Jordan algebra of rank $n$ with the 
symmetric cone $V_+$. Let $Q$ be a permutation invariant proper polyhedral cone in $\Rn$ that is reducible. Let $Q=A(\Rnp)$ 
 as in Theorem \ref{theorem 1}. 
Then,  the spectral cone $\lambda^{-1}(Q)$ is isomorphic to $V_+$. In fact, 
$$\lambda^{-1}(Q)=L_{(a,b)}(V_+),$$
where 
$L_{(a,b)}(x):=(a-b)x+b\,tr(x)e$ is an isomorphism of $V$.
}
\end{proposition}

\begin{proof}
Clearly, $L_{(a,b)}$ is linear on $V$. If $L_{(a,b)}(x)=0$ for some $x$, then, $(a-b)x+b\,tr(x)e=0$. 
Upon taking the trace and noting that $tr(e)=n$, we get $(a-b)tr(x)+b\,tr(x)n=0$. Since $a\neq b(1-n)$, we get
$tr(x)=0$ and hence, $(a-b)x=0$. As $a\neq b$, we have $x=0$. This proves that $L_{(a,b)}$ is an isomorphism of $V$.
Hence, $L_{(a,b)}(V_+)$ is isomorphic to $V_+$. We now show
 that 
$$\lambda^{-1}(Q)=\{(a-b)x+b\,tr(x)e: x\in V_+\}.$$
Let $y\in \lambda^{-1}(Q)$ with its spectral representation $y=\sum y_if_i$, where $y_1,y_2,\ldots, y_n$ are the eigenvalues of $y$ (written in the decreasing order) and $\{f_1,f_2,\ldots,f_n\}$ is a Jordan frame. 
Then, $\lambda(y)=[y_1,y_2,\ldots, y_n]^T=Aq=[(a-b)I+bE]q$ for some $q\geq 0$ in $\Rn$. We see that $y_i=(a-b)q_i+b\tr(q)$ for all $i$ and so $y=(a-b)x+b\tr(x)e$, where $x:=\sum q_if_i\in V_+$. 
Now to see the reverse inclusion, let $y=(a-b)x+b\,tr(x)e$ for some $x\in V_+$. Then, starting with the spectral 
decomposition $x=\sum x_ie_i$ (where $x_1,x_2,\ldots, x_n$ are the eigenvalues of $x$ and $\{e_1,e_2,\ldots,e_n\}$ is a Jordan frame), we see that $y$ has the spectral decomposition $y=\sum \left [ (a-b)x_i+b\tr(x)\right ]e_i$. Letting $y_i:=(a-b)x_i+b\tr(x)$, we see that $y_1,y_2,\ldots, y_n$ are the eigenvalues of $y$. Writing $q=[x_1,x_2,\ldots, x_n]^T$, we see that 
$q\in \Rnp$ (as $x\in V_+$) and $\lambda(y)=(Aq)^{\downarrow}\in A(\Rnp)=Q$.
This completes the proof.
\end{proof}


\end{document}